\documentclass[12pt,fleqn]{article}

\usepackage{geometry}

\usepackage[cp1251]{inputenc}
\usepackage[T2A]{fontenc}
\usepackage[english]{babel}

\usepackage{amsmath}
\usepackage{amssymb}
\usepackage{amsthm}
\usepackage{array}
\usepackage[auth-sc]{authblk}
\usepackage{bm}
\usepackage{cite}
\usepackage{color}
\usepackage{dsfont}
\usepackage{enumerate}
\usepackage{enumitem}
\usepackage{epsf}
\usepackage{euscript}
\usepackage{graphicx} 
\usepackage{indentfirst}
\usepackage{latexsym}
\usepackage{mathrsfs}
\usepackage{mathtools}
\usepackage{tikz}
\usepackage{pgf}

\usetikzlibrary{calc}
\usetikzlibrary{intersections}
\usetikzlibrary{fadings}

\definecolor{cg4}{gray}{0.99}          
\definecolor{cg6}{gray}{0.30}          
\definecolor{cg1}{gray}{0.60}          
\definecolor{cg2}{gray}{0.90}          
\definecolor{cg3}{gray}{0.75}          
\definecolor{cg5}{gray}{0.45}          

\colorlet{c1}{cg1}
\colorlet{c2}{cg2}
\colorlet{c3}{cg3}
\colorlet{c4}{cg4}
\colorlet{c5}{cg5}
\colorlet{c6}{cg6}

\newcommand{\logic}[1]{\mathbf{#1}}
\newcommand{\otuple}[1]{\langle{#1}\rangle}
\newcommand{\kframe}[1]{\mathfrak{#1}}
\newcommand{\kFrame}[1]{\boldsymbol{\kframe{#1}}}
\newcommand{\kmodel}[1]{\mathfrak{#1}}
\newcommand{\kModel}[1]{\boldsymbol{\kmodel{#1}}}
\newcommand{\leftsquare}[1]{%
\begin{tikzpicture}[scale=#1]
\draw (0,0)--(0,1)--(1,1)--(1,0)--cycle;
\draw [fill, color=black!25] (0,0)--(0.5,0.5)--(0,1);
\draw (0,0)--(0.5,0.5)--(0,1);
\draw (0,1)--(0,0);
\end{tikzpicture}
}
\newcommand{\rightsquare}[1]{%
\begin{tikzpicture}[scale=#1]
\draw (0,0)--(0,1)--(1,1)--(1,0)--cycle;
\draw [fill, color=black!25] (1,0)--(0.5,0.5)--(1,1);
\draw (1,0)--(0.5,0.5)--(1,1);
\draw (1,1)--(1,0);
\end{tikzpicture}
}
\newcommand{\upsquare}[1]{%
\begin{tikzpicture}[scale=#1]
\draw (0,0)--(0,1)--(1,1)--(1,0)--cycle;
\draw [fill, color=black!25] (0,1)--(0.5,0.5)--(1,1);
\draw (0,1)--(0.5,0.5)--(1,1);
\draw (1,1)--(0,1);
\end{tikzpicture}
}
\newcommand{\downsquare}[1]{%
\begin{tikzpicture}[scale=#1]
\draw (0,0)--(0,1)--(1,1)--(1,0)--cycle;
\draw [fill, color=black!25] (0,0)--(0.5,0.5)--(1,0);
\draw (0,0)--(0.5,0.5)--(1,0);
\draw (1,0)--(0,0);
\end{tikzpicture}
}
\newcommand{\leftsq} {\mathop{\leftsquare {0.21}}}
\newcommand{\rightsq}{\mathop{\rightsquare{0.21}}}
\newcommand{\upsq}   {\mathop{\upsquare   {0.21}}}
\newcommand{\downsq} {\mathop{\downsquare {0.21}}}
\newcommand{\sleftsq} {\mathop{\leftsquare {0.16}}}
\newcommand{\srightsq}{\mathop{\rightsquare{0.16}}}
\newcommand{\supsq}   {\mathop{\upsquare   {0.16}}}
\newcommand{\sdownsq} {\mathop{\downsquare {0.16}}}

\newcommand{\drawtileflat} [9]
{
\draw [white, opacity = 0, name path = diag 1] #1--#3;
\draw [white, opacity = 0, name path = diag 2] #2--#4;
\draw [name intersections = {of = diag 1 and diag 2, by = {tcenter}}];

\foreach \c in {0.84}
{
  \coordinate (ttl)     at ($(tcenter)+\c*#4-\c*(tcenter)$);
  \coordinate (tbr)     at ($(tcenter)+\c*#2-\c*(tcenter)$);
  \coordinate (tbl)     at ($(tcenter)+\c*#1-\c*(tcenter)$);
  \coordinate (ttr)     at ($(tcenter)+\c*#3-\c*(tcenter)$);
  \coordinate (sttl)    at ($(tcenter)+0.4*(ttl)-0.4*(tcenter)$);
  \coordinate (stbr)    at ($(tcenter)+0.4*(tbr)-0.4*(tcenter)$);
  \coordinate (stbl)    at ($(tcenter)+0.4*(tbl)-0.4*(tcenter)$);
  \coordinate (sttr)    at ($(tcenter)+0.4*(ttr)-0.4*(tcenter)$);
}

\fill [#5, fill opacity=0.75] (ttl)--(tbl)--(stbl)--(sttl)--cycle;
\fill [#6, fill opacity=0.75] (tbl)--(tbr)--(stbr)--(stbl)--cycle;
\fill [#7, fill opacity=0.75] (tbr)--(ttr)--(sttr)--(stbr)--cycle;
\fill [#8, fill opacity=0.75] (ttl)--(ttr)--(sttr)--(sttl)--cycle;

\draw (ttl)--(tbl)--(tbr)--(ttr)--cycle;
\draw (sttl)--(stbl)--(stbr)--(sttr)--cycle;
\draw (ttl)--(sttl);
\draw (ttr)--(sttr);
\draw (tbl)--(stbl);
\draw (tbr)--(stbr);

\node [] at (tcenter) {#9};
}

\newcommand{\Rem}[1]{}

\begin{document}


\title{Tiling problems and complexity of logics\\ (extended version)\thanks{The work on the paper was partially supported by the HSE Academic Fund Programme, Project~\mbox{23-00-022}. The paper is an extended version of the abstracts submitted to the conference SCAN~2023~\cite{RS:SCAN:2023}.}}

\author[1]{Mikhail Rybakov}
\author[2]{Darya Serova}
\affil[1]{IITP RAS, HSE University, Tver State University}
\affil[2]{Tver State University}

\date{}

\maketitle

\thispagestyle{empty}
\newcommand{\insText}[1]{{\color{blue}#1}}
\newcommand{\mText}[1]{{\color{green!60!black}#1}}
\newcommand{\delText}[1]{{\color{red}#1}}

\begin{abstract}
We apply domino problems to give short proofs for some known theorems for the classical predicate logic and to obtain lower bounds for complexity of modal predicate logics defined by Noetherian orders as Kripke frames.
\end{abstract}

\section{Introduction}

Domino, or tiling, problems~\cite{Berger66,Harel86} provide us with a rich tool allowing to
estimate bounds for computational complexity of problems arising in
different fields of mathematics, in particular, in algebra~\cite{CL-algebra-1990,KP-algebra-1999} and
mathematical logic~\cite{BGG97,GOR-1999,RZ01,KKZ05,GKWZ}. Sometimes, properties of tilings of some kind can
be quite easily expressed in a formal language, and their description
can be more elegant than, say, of Turing machines (or other
computational models). Indeed, to describe a tiling, we only have to
say that, for every tile, there are appropriate tiles on the top and
on the right, and that moving right-top or top-right we see the same
tile, while for a Turing machine, to describe just a configuration on
some step of computation, we have to describe a head position, a
state, and symbols stored in tape cells.

Here, we consider two tiling problems, known to be, respectively,
$\Pi^0_1$-complete and $\Sigma^1_1$\nobreakdash-complete, and show examples of
their simulation in first-order theories and logics whose langages are
enriched with some extra expressive means~\cite{GOR-1999} but restricted in the number of individual variables, the number of predicate letters, and their arity. 
Also, we give some examples coming from modal predicate logics defined by special classes of Kripke frames.

\section{Tiling problems we consider}

We may think of a tile as a colored $1 \times 1$ square, with a
fixed orientation. 
A~tile type
$t$ consists of a specification of a color for each edge; we write
$\leftsq t$, $\rightsq t$, $\upsq t$, and $\downsq t$ for the colors
of, respectively, the left, the right, the top, and the bottom edges
of the tiles of type~$t$.

Let $T = \{t_0, \ldots, t_{n}\}$ be a set of tile types. Informally, a
{$T$-tiling} is an arrangement of tiles, whose types are in~$T$, on a grid so that the edge colors of the
adjacent tiles match, both horizontally and vertically; see the
picture below.

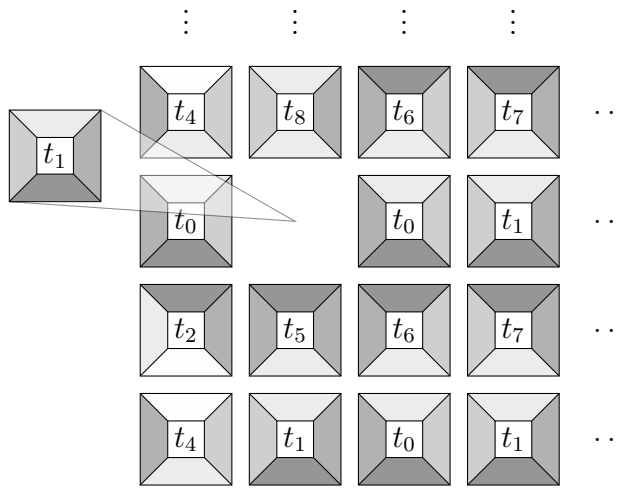
\begin{center}
\begin{tikzpicture}[scale=1.45]

\coordinate (g00)   at (0,0);
\coordinate (g10)   at (1,0);
\coordinate (g20)   at (2,0);
\coordinate (g30)   at (3,0);
\coordinate (g40)   at (4,0);
\coordinate (g50)   at (5,0);
\coordinate (g01)   at (0,1);
\coordinate (g11)   at (1,1);
\coordinate (g21)   at (2,1);
\coordinate (g31)   at (3,1);
\coordinate (g41)   at (4,1);
\coordinate (g51)   at (5,1);
\coordinate (g02)   at (0,2);
\coordinate (g12)   at (1,2);
\coordinate (g22)   at (2,2);
\coordinate (g32)   at (3,2);
\coordinate (g42)   at (4,2);
\coordinate (g52)   at (5,2);
\coordinate (g03)   at (0,3);
\coordinate (g13)   at (1,3);
\coordinate (g23)   at (2,3);
\coordinate (g33)   at (3,3);
\coordinate (g43)   at (4,3);
\coordinate (g53)   at (5,3);
\coordinate (g04)   at (0,4);
\coordinate (g14)   at (1,4);
\coordinate (g24)   at (2,4);
\coordinate (g34)   at (3,4);
\coordinate (g44)   at (4,4);
\coordinate (g54)   at (5,4);
\coordinate (g05)   at (0,5);
\coordinate (g15)   at (1,5);
\coordinate (g25)   at (2,5);
\coordinate (g35)   at (3,5);
\coordinate (g45)   at (4,5);
\coordinate (g55)   at (5,5);

\node [] at ($0.3*(g40)+0.3*(g41)+0.2*(g50)+0.2*(g51)$) {$\ldots$};
\node [] at ($0.3*(g41)+0.3*(g42)+0.2*(g51)+0.2*(g52)$) {$\ldots$};
\node [] at ($0.3*(g42)+0.3*(g43)+0.2*(g52)+0.2*(g53)$) {$\ldots$};
\node [] at ($0.3*(g43)+0.3*(g44)+0.2*(g53)+0.2*(g54)$) {$\ldots$};

\node [] at ($0.3*(g04)+0.3*(g14)+0.2*(g05)+0.2*(g15)$) {$\vdots$};
\node [] at ($0.3*(g14)+0.3*(g24)+0.2*(g15)+0.2*(g25)$) {$\vdots$};
\node [] at ($0.3*(g24)+0.3*(g34)+0.2*(g25)+0.2*(g35)$) {$\vdots$};
\node [] at ($0.3*(g34)+0.3*(g44)+0.2*(g35)+0.2*(g45)$) {$\vdots$};


\drawtileflat{(g00)}{(g10)}{(g11)}{(g01)}{c1}{c2}{c3}{c4}{$t_4$};
\drawtileflat{(g10)}{(g20)}{(g21)}{(g11)}{c3}{c5}{c1}{c2}{$t_1$};
\drawtileflat{(g20)}{(g30)}{(g31)}{(g21)}{c1}{c5}{c3}{c2}{$t_0$};
\drawtileflat{(g30)}{(g40)}{(g41)}{(g31)}{c3}{c5}{c1}{c2}{$t_1$};

\drawtileflat{(g01)}{(g11)}{(g12)}{(g02)}{c2}{c4}{c1}{c5}{$t_2$};
\drawtileflat{(g11)}{(g21)}{(g22)}{(g12)}{c1}{c2}{c1}{c5}{$t_5$};
\drawtileflat{(g21)}{(g31)}{(g32)}{(g22)}{c1}{c2}{c3}{c5}{$t_6$};
\drawtileflat{(g31)}{(g41)}{(g42)}{(g32)}{c3}{c2}{c1}{c5}{$t_7$};

\drawtileflat{(g02)}{(g12)}{(g13)}{(g03)}{c1}{c5}{c3}{c2}{$t_0$};
\drawtileflat{(g22)}{(g32)}{(g33)}{(g23)}{c1}{c5}{c3}{c2}{$t_0$};
\drawtileflat{(g32)}{(g42)}{(g43)}{(g33)}{c3}{c5}{c1}{c2}{$t_1$};

\drawtileflat{(g03)}{(g13)}{(g14)}{(g04)}{c1}{c2}{c3}{c4}{$t_4$};
\drawtileflat{(g13)}{(g23)}{(g24)}{(g14)}{c3}{c2}{c1}{c2}{$t_8$};
\drawtileflat{(g23)}{(g33)}{(g34)}{(g24)}{c1}{c2}{c3}{c5}{$t_6$};
\drawtileflat{(g33)}{(g43)}{(g44)}{(g34)}{c3}{c2}{c1}{c5}{$t_7$};

\drawtileflat{(-1.2,2.6)}{(-0.2,2.6)}{(-0.2,3.6)}{(-1.2,3.6)}{c3}{c5}{c1}{c2}{$t_1$};

\draw [fill=white!10,opacity=0.45] 
       (ttr)--(1.5,2.5)--(tbl)--(tbr)--cycle;

\Rem{

\filldraw [] (g00) circle [radius=1.5pt]   ;
\filldraw [] (g10) circle [radius=1.5pt]   ;
\filldraw [] (g20) circle [radius=1.5pt]   ;
\filldraw [] (g30) circle [radius=1.5pt]   ;
\filldraw [] (g40) circle [radius=1.5pt]   ;

\filldraw [] (g01) circle [radius=1.5pt]   ;
\filldraw [] (g11) circle [radius=1.5pt]   ;
\filldraw [] (g21) circle [radius=1.5pt]   ;
\filldraw [] (g31) circle [radius=1.5pt]   ;
\filldraw [] (g41) circle [radius=1.5pt]   ;

\filldraw [] (g02) circle [radius=1.5pt]   ;
\filldraw [] (g12) circle [radius=1.5pt]   ;
\filldraw [] (g22) circle [radius=1.5pt]   ;
\filldraw [] (g32) circle [radius=1.5pt]   ;
\filldraw [] (g42) circle [radius=1.5pt]   ;

\filldraw [] (g03) circle [radius=1.5pt]   ;
\filldraw [] (g13) circle [radius=1.5pt]   ;
\filldraw [] (g23) circle [radius=1.5pt]   ;
\filldraw [] (g33) circle [radius=1.5pt]   ;
\filldraw [] (g43) circle [radius=1.5pt]   ;

\filldraw [] (g04) circle [radius=1.5pt]   ;
\filldraw [] (g14) circle [radius=1.5pt]   ;
\filldraw [] (g24) circle [radius=1.5pt]   ;
\filldraw [] (g34) circle [radius=1.5pt]   ;
\filldraw [] (g44) circle [radius=1.5pt]   ;
}

\end{tikzpicture}
\end{center}


The fist tiling problem we consider is the following: given a set $T = \{t_0, \ldots, t_{n}\}$ of tile types, we are to determine whether there exists a $T$-tiling $f\colon \mathds{N} \times \mathds{N} \to T$ such that, for every $i, j \in \mathds{N}$,
\begin{itemize}
\item[$(1)$]\quad\!\!\!\!\!\! $\rightsq f(i,j) = \leftsq f(i+1,j)$;
\item[$(2)$]\quad\!\!\!\!\!\! $\upsq f(i,j) = \downsq f(i,j+1)$.
\end{itemize}
This problem is $\Pi^0_1$-complete~\cite{Berger66}. The second tiling problem we consider can be obtained from the first one by adding an extra requirement
\begin{itemize}
\item[$(3)$]\quad\!\!\!\!\!\! the set $\{ j \in \mathds{N} : f(0, j) = t_0 \}$ is infinite,
\end{itemize}
i.e., claiming that there are infinitely many tiles of type $t_0$ in the leftmost column.
This problem is $\Sigma^1_1$-complete~\cite{Harel86}.

\section{Classical theories}

Assume, for simplicity, a classical first-order language with an
infinite supply of monadic predicate letters $P_0,P_1,P_2,\ldots$ and
two binary predicate letters $H$ and~$V$. The intending meaning of
$P_k(x)$ is ``$x$~is placed with a tile of type~$t_k$''; also,
$H(x,y)$ means ``$y$~is to the right of~$x$'', and $V(x,y)$ means
``$y$~is above~$x$''. To describe an $\mathds{N} \times \mathds{N}$
grid, it is sufficient to say
$$
\begin{array}{ccc}
\forall x\exists y\, H(x,y), 
  & \forall x\exists y\, V(x,y), 
  & \forall x\forall y\,(\exists z\, (H(x,z)\wedge V(z,y)) \leftrightarrow \exists z\, (V(x,z)\wedge H(z,y)). 
\end{array}
$$
\pagebreak[3]
Then, we can say that we are given a $T$-tiling:  
\bigskip\nopagebreak[3]

$
\begin{array}{ll}
\mbox{$\bullet$~~Each tile-holder holds a unique tile:}
&
\displaystyle
\forall x\,\bigvee\limits_{\mathclap{i=0}}^n(P_i(x)\wedge \bigwedge\limits_{\mathclap{j\ne i}}\neg P_j(x)).
  \\
\mbox{$\bullet$~~The condition (1) 
is satisfied:}
&
  \displaystyle
  \forall x\,\bigwedge\limits_{\mathclap{i=0}}^n(P_i(x)\to \forall y\,(H(x,y)\to \bigvee\limits_{\mathclap{\srightsq t_i = \hspace{1pt} \sleftsq t_j}}P_j(y))).
  \\
\mbox{$\bullet$~~The condition (2) 
is satisfied:}
&
  \displaystyle
  \forall x\,\bigwedge\limits_{\mathclap{i=0}}^n(P_i(x)\to \forall y\,(V(x,y)\to \bigvee\limits_{\mathclap{\supsq t_i = \hspace{1pt} \sdownsq t_j}}P_j(y))).
\end{array}
$
\bigskip

\noindent%
It is not hard to see that the conjunction of the above formulas is
satisfiable if, and only if, there exists a $T$-tiling
$f\colon \mathds{N} \times \mathds{N} \to T$ satisfying conditions~(1)
and~(2). As a result, the Church's theorem~\cite{Church36} for the
classical first-order logic follows. Since we can simulate all the
predicate letters with a single binary one without adding extra
individual variables~\cite{MR:2022DM,MR:2023LI}, this gives us a short
proof of the known refinement~\cite{TG87} of the Church's theorem: the
satisfiability problem is undecidable for languages with a single
binary predicate letter and three individual variables. Moreover, we
readily obtain undecidability ($\Sigma^0_1$-hardness) for infinite classes of
theories of a binary predicate, again, with three individual variables~\mbox{\cite{MR:2022DM,MR:2023LI}}.

Observe that, with the use of Compactness theorem, the existence of a $T$-tiling satisfying (1) and~(2) is
equivalent to the existance, for every $n\in\mathds{N}$, of an
$n\times n$ tiling with $T$-tiles satisfying (1) and (2) for all
appropriate $i$ and~$j$. Therefore, we can use only finitely many
tile-holders (but their number must be unbounded). This observation
allows us to simulate $T$-tilings on finite models and, thus, to
obtain the Trakhtenbrot's theorem~\cite{Trakhtenbrot50,Trakhtenbrot53}
for satisfiability over finite models. Again, modulo some linguistic
machinations, we obtain undecidability ($\Pi^0_1$-harness) for large
classes of theories of a binary predicate defined by infinite classes
of finite models~\mbox{\cite{MR:2022DM,MR:2023LI}}.

Notice that undecidability of some the theories~--- both $\Sigma^0_1$-hardness and $\Pi^0_1$-hardness~--- follow also from proofs like in~\cite{ELTT:1965,Nies:1996} by means of a general technique described in~\cite{Sper:2016}.

\section{Classical theories with extra non-elementary expressive means}

Having enriched the language with equality and the operator of
transitive closure, we can use the transitive closure $V^+$ of $V$
allowing us to express~$(3)$: 
$$
\exists x\forall y\,(V^+(x,y)\to \exists z\,(z\ne y\wedge V^+(y,z)\wedge P_0(z))).
$$
Notice that equality can be eliminated if we add the condition of
irreflexivity, i.e., 
$$
\forall x\,\neg V(x,x);
$$ 
also, variable $z$ can be
replaced with~$x$. Then, adding the operator of composition $\circ$ of
binary relations, we are able to express that moving right-top and
top-right, we see the same tile, using the formula
$$
\forall x\forall y\,([V\,{\circ}\,H](x,y)\leftrightarrow [H\,{\circ}\,V](x,y)),
$$ 
which contains only two individual variables. Again, using
additional techniques, we can prove that the satisfiability for
languages with a single binary relation, equality, the operators of
transitive closure and composition is $\Sigma^1_1$-hard even for
formulas with two variables~\cite{MR:2022DM}. Sometimes, the operator
of transitive closure can be replaced with the operator asserting
the transitivity of a binary relation~\cite{MR:2023LI}.

\section{Modal predicate logics of Noetherian orders}

The idea to consider together a binary relation and its transitive
closure can be applied also to investigating computational complexity
of non-classical predicate logics using domino problems. We show how
to apply it to modal predicate logics of Kripke
frames that are Noetherian orders.

We shall use the Kripke semantics. To evaluate modal predicate
formulas, we endow a Kripke frame $\kframe{F}=\otuple{W,R}$ with a
system $D=(D_w)_{w\in W}$ of non-empty domains such that
$D_w\subseteq D_v$ whenever $wRv$. To define a model on the resultant
augmented frame $\kFrame{F}=\otuple{\kframe{F},D}$, we endow
$\kFrame{F}$ with a system $I=(I_w)_{w\in W}$ of interpretations of
predicate letters. Then, each world $w$ of a model
$\kModel{M}=\otuple{\kFrame{F},I}$ can be understood as a classical
model $M_w=\otuple{D_w,I_w}$. The truth relation is defined in a
natural way~\cite{GShS}.

Suppose, for simplicity, that $\kframe{F}=\otuple{W,R}$ is a Kripke
frame isomorphic to the structure $1+\omega^\ast$, i.e., there is a
world $w\in W$ seeing an infinite descending $R$-chain isomorphic to
$\omega^\ast$ with the natural strict order. Again, for simplicity,
suppose that all the domains $(D_u)_{u\in W}$ are the same and that
equality $\approx$ is understood as identity. Let $P$ be a binary
predicate letter, corresponding to a serial strict linear order (we
omit transitivity; as we shall show, the transitivity of $P$ follows
from other conditions below):
$$
\begin{array}{ccc}
\forall x\exists y\,P(x,y),
  & \forall x\,\neg P(x,x), 
  & \forall x\forall y\,(P(x,y)\vee x\approx y \vee P(y,x)). 
\end{array}
$$
Let also $V(x,y)$ mean that $y$ is an immediate successor of $x$:
$$
\begin{array}{lcl}
V(x,y) & = & P(x,y)\wedge \neg\exists z\,(P(x,z)\wedge P(z,y)).
\end{array}
$$

For every $x$, we claim that $w$ sees a world ``labelled'' with $x$
using a monadic letter~$L$: $\forall x\,\Diamond L(x)$. Then, we make
a connection between $R^{-1}$ on $W\setminus\{w\}$ and (a predicate
corresponding to) $P$ on $D_w$ via ``labels'':
$$
\forall x\forall y\,(P(x,y) \leftrightarrow \Box(L(y)\to \Diamond L(x))).
$$ 
Notice that the transitivity of $P$ follows. The connection
guarantees that $D_w$ with a relation corresponding to $P$ is
isomorphic ot $\omega$ with the natural strict order. This enables us
to redefine~$V$ using only two variables:
$$
\begin{array}{lcl}
V(x,y) & = & P(x,y)\wedge \Box(L(y)\leftrightarrow \Diamond L(x)\wedge \neg\Diamond\Diamond L(x)). 
\end{array}
$$
Then, we can easily simulate $H$ by moving through all the worlds
except $w$ step-by-step along $R^{-1}$ and choosing the same element
from the domain.

Notice that $P$ is the transitive closure of $V$; therefore, we can
simulate the second tiling problem.
Prior to that, we assert
$$
\begin{array}{lcl}
\forall x\forall y\,(P(x,y)\to \Box P(x,y)) 
  & \mbox{and}
  & \forall x\forall y\,(\neg P(x,y)\to \Box \neg P(x,y)).
\end{array}
$$ 
The formulas are: 
\bigskip

$
\begin{array}{ll}
\mbox{$\bullet$~~Uniqueness:} & \displaystyle
    \forall x\forall y\,\Box \bigvee_{\mathclap{k=0}}^n(P_k(y)\wedge \bigwedge\limits_{\mathclap{j\ne k}}\neg P_{j}(y)).
    \smallskip\\
\mbox{$\bullet$~~Condition (1):} & \displaystyle
    \forall x\forall y\,\Box\bigwedge\limits_{\mathclap{k=0}}^n(L(y)\wedge P_k(x)\to \Box(\exists x\,(V(x,y)\wedge L(x))\to \bigvee\limits_{\mathclap{\mathop{\srightsq} t_j=\hspace{1pt}\mathop{\sleftsq} t_k}} P_{j}(x))).
    \smallskip\\
\mbox{$\bullet$~~Condition (2):} & \displaystyle
    \forall x\forall y\,(V(x,y)\to \Box\bigwedge\limits_{\mathclap{k=0}}^n(P_{k}(x)\to
    \bigvee\limits_{\mathclap{\mathop{\supsq} t_k=\hspace{1pt}\mathop{\sdownsq} t_j}}P_{j}(y))).
    \smallskip\\
\mbox{$\bullet$~~Condition (3):} & \displaystyle
    \forall x\exists y\,(P(x,y) \wedge \Box(\Box\bot\to P_{0}(y))).
    \phantom{\bigwedge\limits^n}
\end{array}
$
\bigskip

Some further linguistic machinations enable us to simulate the tiling
problem in a class of Noetherian orders (strict or non-strict) lying
between the class of all such orders and the class of all linear ones,
using two monadic predicate letters, a proposition letter, and two
individual variables (with some variations)~\cite{MR:2022IGPL}. As a
result, the modal predicate logic of every such class is shown to be
$\Pi^1_1$-hard in such languages. Then, we readily obtain Kripke
incompleteness of all logics lying between $\logic{QwGrz}$ and
$\logic{QwGL.3.bf}$ or between $\logic{QwGrz}$ and
$\logic{QGrz.3.bf}$~\cite{MR:2022IGPL}. This, in particular, gives us
another proof of Kripke incompleteness
of~$\logic{QGL}$~\cite{Montagna84}, even in quite poor languages.


\section{Some remarks and further results}

Examples of the use of tiling problems for obtaining results on the
algorithmic complexity of various logics, both propositional and
predicate, can be found
in~\cite{BGG97,RZ01,GKWZ,KKZ05,RSh20AiML,RShJLC21c,RSh:LI:2021,MR:2022IGPL,MR:2023TSU}. In
particular, the tiling problems considered here can be used to obtain
complexity results for theories of trees~\cite{MR:2023TSU} and to prove
that modal predicate logics whose Kripke frames are Noetherian orders
are $\Pi^1_1$-hard in rather poor languages~\cite{MR:2022IGPL}.

%





\end{document}